\documentclass[12pt]{article}

\newtheorem{thm}{Theorem}[section]
\newtheorem{prop}[thm]{Proposition}
\newtheorem{lm}[thm]{Lemma}
\newtheorem{cor}[thm]{Corollary}

\newtheorem{de}{Definition}[section]

\usepackage{amsfonts}
\usepackage{amsmath}
\usepackage[american]{babel}

\begin{document}

\begin{title}
{Functional inequalities\\ related to the Rogers-Shephard inequality}
\end{title}


\author{Andrea Colesanti
\footnote{Dipartimento di Matematica ``U. Dini'', 
viale Morgagni 67/A, 50134 Firenze, Italy;
colesant@math.unifi.it}} 
\date{}
\maketitle

\begin{center}
{\em Dedicated to Rolf Schneider in occasion of his 65th birthday}
\end{center}

\bigskip

\begin{abstract} 
\noindent For a real-valued nonnegative and log-concave
  function $f$ defined in ${\bf R}^n$, we introduce a notion of
  {\em difference function} $\Delta f$; the difference function
  represents a functional analog on the {\em difference body} $K+(-K)$ of a
  convex body $K$. We
  prove a sharp inequality which bounds the integral of $\Delta f$ from above, in
  terms of the integral of $f$ and we characterize equality conditions. The
  investigation is extended to an analogous notion of difference function for
  $\alpha$-concave functions, with $\alpha<0$. In this case also, we prove an
  upper bound for the integral of the $\alpha$-difference function of $f$ in
  terms of the integral of $f$; the bound is proved to be sharp in the case
  $\alpha=-\infty$ and in the one dimensional case. 
\end{abstract}

\section{Introduction and results}

Convex geometry has among its most important achievements several
well-known inequalities: the Brunn-Minkowski inequality, the
Aleksandrov-Fenchel inequality, the Blaschke-Santal\'o inequality and
many others. 

A recent development in this field, which involves other areas of
mathematics and notably functional analysis, consists in the
interpretation of results having an unmistakable geometric nature, by
an analytic, or better functional, point of view. 

Some aspects of this interaction between convex geometry and analysis
are very well described in the survey paper \cite{Gardner} devoted to
the Brunn-Minkowski inequality. In particular, \S 7 of this paper presents the
Pr\'ekopa-Leindler inequality which has to be considered
the functional counterpart of the Brunn-Minkowski inequality. 
This is probably the most enlightening example of the phenomenon of
'translation' from geometry to analysis that we have mentioned
before; let us briefly see how this translation is made. The
Brunn-Minkowski inequality, in its most general form, claims that if
$A_0$, $A_1$ and $A$ are measurable subsets of ${\bf R}^n$ such that for
some $t\in[0,1]$ 
\begin{equation*}
A\supset A_t:=(1-t)A_0+tA_1=\{(1-t)x+ty\,|\,x\in A_0\,,\,y\in A_1\}\,,
\end{equation*}
then
\begin{equation}
\label{I.0d}
[V_n(A)]^{1/n}\ge(1-t)[V_n(A_0)]^{1/n}+t[V_n(A_1)]^{1/n}\,,
\end{equation}
where $V_n(\cdot)$ denotes the $n$-dimensional volume (i.e. the Lebesgue
measure in ${\bf R}^n$). An equivalent form of
(\ref{I.0d}) is 
\begin{equation}
\label{I.0f}
V_n(A)\ge[V_n(A_0)]^{1-t}[V_n(A_1)]^t\,,
\end{equation}
for every $A_0$ and $A_1$, $t\in[0,1]$ and $A\supset A_t$; this is also
referred to as the {\em multiplicative form} of the Brunn-Minkowski inequality (a proof of
the equivalence between (\ref{I.0d}) and (\ref{I.0f}) can be found
in \cite{Gardner}, \S 7). 

The Pr\'ekopa-Leindler inequality states that if $f_0$, $f_1$ and $f$
are measurable and nonnegative functions defined in ${\bf R}^n$, such
that for some $t\in[0,1]$
\begin{equation}
\label{I.0g}
f(z)\ge f_t(z):=
\sup\{f_0(x)^{1-t}f_1(y)^t\,|\,(1-t)x+ty=z\}
\,,\quad\forall\,z\in{\bf R}^n\,,
\end{equation}
then
\begin{equation}
\label{I.0h}
\int_{{\bf R}^n}f(z)\,dz\ge
\left(\int_{{\bf R}^n}f_0(x)\,dx\right)^{1-t}\,
\left(\int_{{\bf R}^n}f_1(y)\,dy\right)^t\,.
\end{equation}
Notice that if $f_0$ and $f_1$ are characteristic functions of measurable sets 
$A_0$ and $A_1$ respectively, then $f_t$ is the characteristic
function of $K_t$, so that (\ref{I.0h}) implies (\ref{I.0f}).

In the passage from Brunn-Minkowski inequality to
Pr\'ekopa-Leindler inequality, the convex linear combination $A_t$ of
sets $A_0$ and $A_1$ is replaced by the interpolation $f_t$ of
$f_0$ and $f_1$ defined in (\ref{I.0g}). An equivalent
definition of $f_t$ is
\begin{equation*}
f_t(z)=e^{-v_t(z)}
\end{equation*}
where
\begin{equation*}
v_t(z):=
\inf\{(1-t)v_0(x)+tv_1(y)\,|\,(1-t)x+ty=z\}
\end{equation*}
and $v_0=-\log f_0$, $v_1=-\log f_1$. The function $v_t$ is the {\em
  infimal convolution} of $v_0$ and $v_1$, introduced by Rockafellar (in
  the case of convex functions) in \cite{Rockafellar}.

\bigskip

Another result that goes in the same direction is the functional
form of the Blaschke-Santal\'o inequality, proved by
Ball in \cite{Ball} and recently extended by Artstein, Klartag and Milman
in \cite{Artstein-Klartag-Milman}. We are not going to describe this
result in details, but we want to underline a difference with respect to
the previous example. The Brunn-Minkowski and the Pr\'ekopa-Leindler
inequality are valid for measurable sets and nonnegative measurable
functions respectively. On the contrary, the validity of Blaschke-Santal\'o
inequality is restricted to convex bodies and, correspondingly, its
functional version holds in the class of nonnegative {\em log-concave}
functions. Log-concavity appears to be a natural adaptation to
the functional setting of the notion of convexity for sets (notice that
if $f$ is the characteristic function of a set $A$, then $A$ is convex
if and only if $f$ is log-concave).  

\bigskip

In this paper we deal with inequalities which can be viewed
as functional forms of a third inequality of
convex geometry, the Rogers-Shephard inequality (see
\cite{Rogers-Shephard1} and \S 7.3 in \cite{Schneider}). This inequality
provides an optimal upper bound for the volume of the difference body of
a convex body $K$, in terms of the 
volume of $K$. Let us recall that an $n$-dimensional convex body $K$ ($n\ge
2$) is a compact convex subset of ${\bf R}^n$; the difference body $DK$
of $K$ is defined as
\begin{equation*}
\label{I.0a}
DK=K+(-K)=\{x+y\,|\,x\in K\,,\,-y\in K\}\,;
\end{equation*}
$DK$ is also a convex body and it is symmetric with respect to the origin. The
Rogers-Shephard inequality states that for every $n$-dimensional convex body $K$
\begin{equation}
\label{I.0b}
V_n(DK)\le\binom{2n}{n}V_n(K)\,.
\end{equation}
This inequality is optimal, indeed it
becomes an equality when $K$ is a simplex (and only in this case).
$V_n(DK)$ can be also estimated from below, using the Brunn-Minkowski
inequality:
\begin{equation}
\label{I.0c}
V_n(DK)\ge2^nV_n(K)\,.
\end{equation}
The validity of Rogers-Shephard inequality is restricted to
convex bodies, as simple examples (even in dimension one) show. 

\bigskip

In the sequel, we introduce the notion of {\em difference
function} of a nonnegative function $f$ defined in
${\bf R}^n$, based on the functional
interpolation that we have seen in (\ref{I.0g}):
roughly speaking, the difference function of $f$ is the interpolation of
$f(x)$ and $f(-x)$, with $t=\frac{1}{2}$.

\begin{de} Let $f$ be a real-valued, nonnegative function, defined in
  ${\bf R}^n$. The difference function $\Delta f$ of $f$ is defined as
\begin{equation}
\label{I.0m}
\Delta f(z)=\sup\left\{\sqrt{f(x)f(-y)}\;|\;x,y\in{\bf R}^n\,,\,\frac{1}{2}(x+y)=z
\right\}\,,\quad\forall\,z\in{\bf R}^n\,.
\end{equation}
\end{de}
$\Delta f$ is an even function.
If we use Pr\'ekopa-Leindler inequality we obtain 
\begin{equation*}
\int_{{\bf R}^n}f(x)\,dx\le
\int_{{\bf R}^n}\Delta f(z)\,dz\,.
\end{equation*}
In other words, Pr\'ekopa-Leindler provides a lower bound for
the integral of $\Delta f$ as well as Brunn-Minkowski inequality
provides a lower bound for the volume of the difference
body. Our purpose is to prove a
corresponding upper bound for the integral of $\Delta f$ in terms of the
integral of $f$ provided that $f$ is log-concave, which
corresponds to the convexity assumption in the Rogers-Shephard
inequality. 

\begin{de}
Let $f$ be a real-valued nonnegative function defined in ${\bf R}^n$; we
say that $f$ is log-concave in ${\bf R}^n$ if
\begin{equation*}
\label{I.1}
f((1-t)x+ty)\ge f(x)^{1-t}f(y)^t\,,\quad
\forall x,y\in{\bf R}^n\,,\quad\forall t\in[0,1]\,.
\end{equation*}
\end{de}
This is equivalent to say that the function
\begin{equation*}
\label{I.2}
v=-\log(f)\,:{\bf R}^n\longrightarrow(-\infty,+\infty]
\end{equation*}
is convex in ${\bf R}^n$.

If $f$ is log-concave, then its difference function is also log-concave;
indeed, as we have seen before, 
\begin{equation*}
\Delta f(z)=e^{-\delta v(z)}
\end{equation*}
where 
\begin{equation}
\label{I.0aa}
\delta v(z)=\inf\left\{\frac{v(x)+v(-y)}{2}\,|\,x,y\in{\bf R}^n\,,\,
\frac{1}{2}(x+y)=z\right\}
\end{equation}
(here $v=-\log(f)$ as above). By the convexity of $v$, $\delta v$ is convex 
(see \cite{Rockafellar}, Chapter 5); consequently, $\Delta f$ is log-concave.

We are ready to state our main results.

\begin{thm} 
\label{teo.I.1}
Let $f$ be a real-valued, nonnegative and log-concave function
defined in ${\bf R}^n$. Then
\begin{equation}
\label{I.4}
\int_{{\bf R}^n}\Delta f(z)\,dz\le
2^n\,\int_{{\bf R}^n}f(x)\,dx\,.
\end{equation}
\end{thm}

Inequality (\ref{I.4}) can not be improved as it is showed by the following result.

\begin{thm} 
\label{teo.I.2}
Let $g$ be defined as follows
\begin{equation*}
\label{I.5}
g(x)=g(x_1,\dots,x_n)=
\left\{
\begin{array}{ll}
e^{-(x_1+\dots+x_n)}\;&\mbox{if}\;x_i\ge 0\quad\forall i=1,\dots,n\,,\\
0&\mbox{otherwise.}
\end{array}
\right.
\end{equation*}
Then 
\begin{equation*}
\label{I.6}
\Delta g(z)=\Delta g(z_1,\dots,z_n)=e^{-(|z_1|+\dots+|z_n|)}\,,\quad\forall z\in{\bf R}^n
\end{equation*}
and 
\begin{equation*}
\label{I.7}
\int_{{\bf R}^n}\Delta g(z)\,dz=2^n\int_{{\bf R}^n}g(x)\,dx\,.
\end{equation*}
\end{thm}

\bigskip

Let $g$ be as in the previous theorem; if $A$ is a non-singular square
matrix of order $n$, 
$x_0\in{\bf R}^n$ and $C>0$, then the function
\begin{equation}
\label{I.4b}
h(x)=Cg(Ax+x_0)\,,\quad x\in{\bf R}^n
\end{equation}
is also an extremal function for inequality (\ref{I.4}), i.e.
\begin{equation*}
\int_{{\bf R}^n}\Delta h(x)\,dx=2^n\int_{{\bf R}^n}h(x)\,dx\,.
\end{equation*}
The next result states that functions of the form (\ref{I.4b}) exhaust
the family of extremal functions for inequality (\ref{I.4}). 
Let $f$ be a real-valued, nonnegative and log-concave function defined in
${\bf R}^n$. We set
\begin{equation*}
P_f=\{x\in{\bf R}^n\,|\,f(x)>0\}\,.
\end{equation*}
By the log-concavity of $f$, the set $P_f$ is convex; we denote by 
${\rm int}(P_f)$ and ${\rm cl}(P_f)$ its interior and its closure 
respectively. 

\begin{thm}\label{teo.II.1}
Let $f$ be a real-valued, nonnegative and log-concave function.
Assume that $f\in L^1({\bf R}^n)$, ${\rm int}(P_f)\ne\emptyset$ and   
\begin{equation*}
\int_{{\bf R}^n}\Delta f(z)\,dz=
2^n \int_{{\bf R}^n}f(x)\,dx\,.
\end{equation*}
Then there exists a non-singular matrix $A$ of order $n$, $C>0$ and
$x_0\in{\bf R}^n$ such that if $h$ is defined by
\begin{equation*}
h(x)=Cf(Ax+x_0)\,,
\end{equation*}
then:\\
i) ${\rm int}(P_h)=\{x=(x_1,\dots,x_n)\,:\,x_i>0\,,\,i=1,\dots,n\}$;\\
ii) $h(z)=h(z_1,\dots,,z_n)=e^{-(z_1+\dots+z_n)}$ in ${\rm int}(P_h)$.
\end{thm}

\bigskip


\bigskip

Theorems \ref{teo.I.1}, \ref{teo.I.2} and \ref{teo.II.1} 
are proved in \S \ref{Paragrafo2}. In the proof of Theorem
\ref{teo.I.1} we use the same idea as in the original proof of the
Rogers-Shephard inequality (see \cite{Rogers-Shephard1}), which renders
the argument rather simple; on the other hand the characterization of
equality conditions (Theorem \ref{teo.II.1}) requires a rather delicate
(and lengthy) argument.  


\bigskip

In \S \ref{Paragrafo3} we deduce from (\ref{I.4}) an inequality which gives an
optimal upper bound for the volume of the convex hull of $K$ and $-K$, where
$K$ is a convex body containing the origin. This result was already known and
its original proof is due to Rogers and Shephard, see
\cite{Rogers-Shephard2}.  
 
\bigskip

In \S \ref{Paragrafo4} and \S\ref{Paragrafo5} we prove some extensions of Theorem \ref{teo.I.1} 
arising from the following consideration. The definition of difference function
(\ref{I.0m}) of a function $f$ is based on {\em geometric means} of values of
$f$; would it be reasonable to consider other means? For $\alpha\in{\bf R}$
and $f\,:\,{\bf R}^n\rightarrow[0,\infty]$, let us 
define the $\alpha$-th difference function of $f$ 
\begin{equation*}
\Delta_\alpha f(z)=\sup\left\{M_\alpha(f(x),f(-y))\,:\,\frac{1}{2}(x+y)=z\right\}\,,\quad z\in{\bf R}^n\,,
\end{equation*}
where
\begin{equation*}
M_\alpha(a,b)=
\left\{
\begin{array}{lll}
&\left[\displaystyle{\frac{a^\alpha+b^\alpha}{2}}\right]^{1/\alpha}\;&\mbox{if $\alpha\ne 0$,}\\
\\
&\sqrt{ab}\;&\mbox{if $\alpha=0$,} 
\end{array}
\right.\quad\forall a,b\ge0
\end{equation*}
(if $\alpha<0$ and $ab=0$, put $M_\alpha(a,b)=0$). The problem now is to prove an inequality of the form
\begin{equation}
\label{aggiuntaI.1}
\int_{{\bf R}^n}\Delta_\alpha f(z)\,dz\le C\int_{{\bf R}^n} f(x)\,dx
\end{equation}
for every $f$ in a suitable class of functions, where $C$ is a
constant independent of $f$, and to determine the best possible constant
$C(n,\alpha)$ for which it is true. Note that in the case $\alpha=0$ the solution is given by
Theorems \ref{teo.I.1} and \ref{teo.I.2}. As we will see, the problem is
meaningful only for $\alpha\le 0$ and for $f$ such that $f^\alpha$ is convex. Under these assumptions we prove inequality
(\ref{aggiuntaI.1}) and we determine the optimal constant in two cases; the
first is $n\ge 1$ and $\alpha=-\infty$ (the $-\infty$ mean is the minimum),
and the second, treated in \S\ref{Paragrafo5}, is $n=1$ and $\alpha\le 0$.

\section{Proof of the main results}
\label{Paragrafo2} 

As a general fact, we notice that a nonnegative log-concave function $f$
is  measurable; indeed, for every $s\ge0$ the set
$\{x\in{\bf R}^n\,:\,f(x)\ge s\}$
is convex, and in particular it is measurable.

\noindent{\bf Proof of Theorem \ref{teo.I.1}.} 
We may assume that $f\in L^1({\bf R}^n)$. Let us define 
the function
\begin{equation*}
F(z)=\int_{{\bf R}^n}f(x-z)f(x)\,dx\,.
\end{equation*}
We fix $z\in{\bf R}^n$. Let $x_i, y_i\in{\bf R}^n$, $i\in{\bf N}$, be two 
sequences such that
\begin{equation*}
z=\frac{1}{2}(x_i+y_i)\,,\;\forall i\in{\bf N}\quad
\mbox{and}\quad
\Delta f(z)=\lim_{i\to\infty}\sqrt{f(x_i)f(-y_i)}\,.
\end{equation*}
As $f$ is log-concave, we have for $x\in{\bf R}^n$ and $i\in{\bf N}$
\begin{equation*}
f(x-z)=f\left(\frac{1}{2}(2x-x_i)+\frac{1}{2}(-y_i)\right)
\ge\sqrt{f(2x-x_i)f(-y_i)}\,,
\end{equation*}
and similarly
\begin{equation*}  
f(x)\ge\sqrt{f(2x-x_i)f(x_i)}\,.
\end{equation*}
Consequently, for every $i\in{\bf N}$,
\begin{eqnarray*}
F(z)\ge\sqrt{f(x_i)f(-y_i)}\int_{{\bf R}^n}f(2x-x_i)\,dx
=\sqrt{f(x_i)f(-y_i)}\;\frac{1}{2^n}\int_{{\bf R}^n}f(x)\,dx\,.
\end{eqnarray*}
Letting $i$ tend to infinity we obtain, for every $z\in{\bf R}^n$,
\begin{equation}
\label{II.0}
F(z)\ge \Delta f(z)\frac{1}{2^n}\int_{{\bf R}^n}f(x)\,dx\,,
\end{equation}
so that
\begin{equation*}
\int_{{\bf R}^n}F(z)\,dz\ge\frac{1}{2^n}
\int_{{\bf R}^n}\Delta f(z)\,dz
\int_{{\bf R}^n}f(x)\,dx\,.
\end{equation*}
On the other hand
\begin{equation*}
\int_{{\bf R}^n}F(z)\,dz=
\left(\int_{{\bf R}^n}f(x)\,dx\right)^2\,,
\end{equation*}
this concludes the proof.

\begin{flushright}$\Box$
\end{flushright}

\bigskip

\noindent{\bf Proof of Theorem \ref{teo.I.2}.} By the definition of
difference function we may write
\begin{equation*}
\Delta g(z)=\sup_{x\in{\bf R}^n}\sqrt{f(x)f(x-2z)}\,.
\end{equation*}
For $z=(z_1,\dots,z_n)\in{\bf R}^n$ we define the set
\begin{equation*}
B_z=\{x=(x_1,\dots,x_n)\in{\bf R}^n\,|\,x_i\ge\max\{0,2z_i\}\;{\rm for}\;i=1,\dots,n\}\,.
\end{equation*}
The expression $g(x)g(x-2z)$ vanishes if $x$ belongs to the complement of 
$B_z$, while for $x\in{B_z}$ its value is
\begin{eqnarray*}
g(x)g(x-2z)&=&
\left[e^{-(x_1+\dots+x_n)}\,e^{-[(x_1-2z_1)+\dots+(x_n-2z_n)]}\right]^{1/2}\\
&=&
e^{z_1+\dots+z_n}\,e^{-(x_1+\dots+x_n)}\,.
\end{eqnarray*}
Thus, for $z\in{\bf R}^n$,
\begin{eqnarray*}
\Delta g(z)&=&e^{z_1+\dots+z_n}\sup_{x\in B_z}e^{-(x_1+\dots+x_n)}
=e^{z_1-2\max\{0,z_1\}+\dots+z_n-2\max\{0,z_n\}}
=e^{-(|z_1|+\dots+|z_n|)}\,.
\end{eqnarray*}

Finally, notice that
\begin{equation*}
\int_{{\bf R}^n}g(x)\,dx=\left(\int_0^\infty
e^{-t}\,dt\right)^n=1\quad
\mbox{and}\quad
\int_{{\bf R}^n}\Delta g(x)\,dx=\left(\int_{-\infty}^\infty e^{-|t|}\,dt\right)^n
=2^n\,.
\end{equation*}

\begin{flushright}
$\Box$
\end{flushright}

\bigskip

For the proof of Theorem \ref{teo.II.1} we need two auxiliary results
stated as lemmas. Let us recall that for a real-valued
nonnegative function $f$ defined in ${\bf R}^n$, $P_f=\{x\,:\,f(x)>0\}$. 

\begin{lm}
\label{lemma.II.1}
Let $f$ be a real-valued, nonnegative and log-concave function defined in 
${\bf R}^n$. Assume that the set $P_f$ has non-empty interior. Then the function
\begin{equation*}
{\bar f}(x)=\left\{
\begin{array}{lll}
f(x)\;&{\rm if}\; x\in{\rm int}(P_f)\\
\limsup_{y\to x,y\in{\rm int}(P_f)}f(y) \;&{\rm if}\; x\in\partial P_f\\
0   \;&{\rm otherwise,}
\end{array}
\right.
\end{equation*}
is log-concave in ${\bf R}^n$.
\end{lm}

\noindent{\bf Proof.} We have to prove that
\begin{equation}
\label{II.1}
{\bar f}((1-t)x+ty)\ge{\bar f}(x)^{1-t}{\bar f}(y)^t\,,
\end{equation}
for every $x,y\in{\bf R}^n$ and for every $t\in[0,1]$. If either $x$ or $y$ 
belongs to ${\bf R}^n\setminus{\rm cl}(P_f)$ the inequality is obvious. So, 
assume that $x,y\in{\rm cl}(P_f)$. Let $x_i$, $y_i$, $i\in{\bf N}$, be two 
sequences such that
\begin{equation*}
x_i,y_i\in{\rm int}(P_f)\;\forall i\in{\bf N}\,,\;
\lim_{i\to\infty}x_i=x\,,\;
\lim_{i\to\infty}y_i=y\,,
\end{equation*}
and 
\begin{equation*}
\lim_{i\to\infty}f(x_i)={\bar f}(x)\,,\;
\lim_{i\to\infty}f(y_i)={\bar f}(y)\,.
\end{equation*}
Then 
\begin{equation*}
{\bar f}((1-t)x+ty)\ge\limsup_{i\to\infty}f((1-t)x_i+ty_i)\ge
\limsup_{i\to\infty}f(x_i)^{1-t}f(y_i)^t
={\bar f}(x)^{1-t}{\bar f}(y)^t\,,
\end{equation*}
where, in the case $(1-t)x+ty\in{\rm int}(P_f)$, the first inequality is
due to the continuity of $f$ in ${\rm int}(P_f)$.

\begin{flushright}
$\Box$
\end{flushright}

\bigskip

\begin{lm} 
\label{lemma.II.2}
Let $f$ be a real-valued, nonnegative and log-concave function defined 
in ${\bf R}^n$. Assume that the set $P_f$ has nonempty interior. Then
\begin{equation*}
\sup_{x\in{\bf R}^n}f(x)=\mathop{{\rm ess\,sup}}_{x\in{\bf R}^n}f(x)\,.
\end{equation*}
\end{lm}

\noindent
{\bf Proof.} Let $v=-\log(f)$; the assertion of the Lemma is equivalent to
\begin{equation*}
\inf_{x\in{\bf R}^n}v(x)=\mathop{\rm ess\,inf}_{x\in{\bf R}^n}v(x)\,.
\end{equation*}
Let $s'$ be such that the set $S=\{x\in{\bf R}^n\,:\,v(x)\le s'\}$ is 
non-empty and its Lebesgue measure is 0; we prove that $v$ is constant  
in $S$. By contradiction, let $x'\in S$ be such that $v(x')<s'$;
$S$ is a convex set of null measure, then there exists an hyperplane $H$
such that $S\subset H$.  
The interior of $P_f$ is non-empty, so that
we may take a point $y\in P_f\setminus H$. The function
\begin{equation*}
\phi(r)=v((1-r)y+rx')\,,\quad r\in[0,1]\,,
\end{equation*}
is convex; on the other hand 
\begin{equation*}
\phi(r)\ge s'\;{\rm for}\; r\in[0,1)\,,\;\phi(1)=f(x')<s'\,;
\end{equation*}
this contradicts the convexity of $\phi$. We have then proved that
\begin{equation*}
v(x)=s'\quad\forall x\in\{x\in{\bf R}^n\;:\;v(x)\le s'\}
\end{equation*}
but this clearly imply that $s'=\inf_{{\bf R}^n}v$. We deduce that for every 
$s>\inf\{v(x)\,|\,x\in{\bf R}^n\}$ the set
$\{x\in{\bf R}^n\;:\;v(x)\le s\}$ has interior points and therefore
positive measure. This concludes the proof of the lemma. 

\begin{flushright}
$\Box$
\end{flushright}

\bigskip

\noindent{\bf Proof of Theorem \ref{teo.II.1}.} The proof is divided into several 
steps.

{\em Step 1.} For convenience we write $P$ instead of $P_f$. For $t>0$
we define
$$
P_t=\{x\in{\bf R}^n\,:\,f(x)\ge t\}\,.
$$
$P_t$ is a convex set. We will prove that:
{\em i}) $P_t$ is bounded for every $t>0$;
{\em ii}) $\sup_{{\bf R}^n}f<\infty$.

{\em Proof of i)}. Since 
$f\in L^1({\bf R}^n)$, the Lebesgue measure of $P_t$ is finite for every
$t>0$. As ${\rm int}(P)\ne\emptyset$, there exists $\bar t$ be such that 
${\rm int}(P_{\bar t})\ne\emptyset$. A convex set with interior points and finite 
measure is bounded, then $P_{\bar t}$ is bounded. By the inclusion 
$P_t\subset P_{\bar t}$, for $t\ge{\bar t}$, the same conclusion holds for 
$t\ge{\bar t}$. On the other hand, for every $t\le{\bar t}$, 
${\rm int}(P_t)\supset{\rm int}(P_{\bar t})\ne\emptyset$, so that by the same 
argument as above, $P_t$ is bounded also for $t\le{\bar t}$.

{\em Proof of ii)}. Let $0<t<\sup_{{\bf R}^n}f$; by the previous step, $P_t$ is 
bounded so that the convex function $v=-\log(f)$ is bounded from below in $P_t$ and
consequently $f$ is bounded from above in $P_t$. The conclusion follows as 
$\sup_{{\bf R}^n}f=\sup_{P_t}f$.

Let us note that a further consequence of {\em i}) is that 
$$
\lim_{|x|\to\infty}f(x)=0\,.
$$

\bigskip

{\em Step 2}. We prove that, without loss of generality, we may assume that
for every $x\in\partial P$
\begin{equation}
\label{II.1b}
f(x)=\limsup_{y\to x,y\in{\rm int}(P)}f(y)\,.
\end{equation}
{\em Proof.} For a nonnegative function $h$ defined in ${\bf R}^n$, let us define
\begin{equation*}
\Delta ^*h(z)=\mathop{\rm ess\,sup}_{x\in{\bf
    R}^n}\sqrt{h(x)h(x-2z)}\,,\;z\in{\bf R}^n\,.
\end{equation*}
In particular, if $h$ is log-concave, for every $z$ the function 
$x\rightarrow\sqrt{h(x)h(x-2z)}$ is also log-concave, thus, by Lemma 
\ref{lemma.II.2}
\begin{equation*}
\Delta ^*h(z)=\Delta h(z)\,,\quad\forall z\in{\bf R}^n\,.
\end{equation*}

We set 
\begin{equation*}
f^*(x)=\left\{
\begin{array}{ll}
f(x)\;&{\rm if}\; x\in{\rm int}(P)\\
0   \;&{\rm otherwise;}
\end{array}
\right.
\end{equation*}
it is easy to verify that $f^*$ is log-concave, so that
\begin{equation}
\label{II.1c}
\Delta^*f^*(z)=\Delta f^*(z)\,,\quad\forall z\in{\bf R}^n\,.
\end{equation}
Moreover, in \cite{Brascamp-Lieb}
(Theorem 2 of the Appendix) it is proved that  
\begin{equation*}
\Delta ^*f^*(z)=\Delta f(z)\,,\quad\forall z\in{\bf R}^n\,.
\end{equation*}
Let us also define  
\begin{equation*}
{\bar f}(x)=\left\{
\begin{array}{lll}
f(x)\;&{\rm if}\; x\in{\rm int}(P)\\
\limsup_{y\to x,y\in{\rm int}(P)}f(y) \;&{\rm if}\; x\in\partial P\\
0   \;&{\rm otherwise.}
\end{array}
\right.
\end{equation*}
By Lemma \ref{lemma.II.1}, $\bar f$ is log-concave and then 
$\Delta ^*{\bar f}=\Delta {\bar f}$. On the other hand
\begin{equation*}
{\bar f}=f^*\;\mbox{a.e. in}\;{\bf R}^n\,,
\end{equation*}
so that for every $z\in{\bf R}^n$
\begin{equation*}
{\bar f}(x){\bar f}(x-2z)=
f^*(x)f^*(x-2z)\;\mbox{ for a.e. $x$ in ${\bf R}^n$.}
\end{equation*}
As a consequence we have
\begin{equation}
\label{II.1d}
\Delta f(z)=\Delta ^*f^*(z)=\Delta ^*{\bar f}(z)=\Delta {\bar f}(z)\,,\quad\forall z\in{\bf R}^n\,.
\end{equation}
We infer that if $f$ verifies the assumptions of Theorem \ref{teo.II.1},
then the same does ${\bar f}$; moreover $P_{\bar f}=P$ and ${\bar f}=f$
in ${\rm int}(P)$, so
that if we prove Theorem \ref{teo.II.1} for $\bar f$, we automatically
prove it for $f$ as well. In the rest of the proof we replace $f$ by
$\bar f$, and for simplicity we continue to write $f$ instead of $\bar f$. 

\bigskip

{\em Step 3.} We prove that for every $z\in{\bf R}^n$ there exist
$x,y\in{\bf R}^n$ such that 
\begin{equation*}
\label{II.2}
\frac{1}{2}(x+y)=z\,,\quad \Delta f(z)=\sqrt{f(x)f(-y)}\,.
\end{equation*}

{\em Proof.} If $\Delta f(z)=0$, this is obvious. Let $\Delta f(z)>0$; using the
notation introduced in the previous step, we have by (\ref{II.1c}) and
(\ref{II.1d}) 
$$
\Delta f(z)=\Delta^*f^*(z)=\Delta f^*(z)\,.
$$
Consequently, we may find two 
sequences $x_i, y_i\in{\bf R}^n$, $i\in{\bf N}$, such that
\begin{equation}
\label{II.3}
x_i,-y_i\in{\rm int}(P)\,,\;
\frac{1}{2}(x_i+y_i)=z\,,\;\forall i\in{\bf N}\,,\quad
\Delta f(z)=\lim_{i\to\infty}\sqrt{f(x_i)f(-y_i)}\,.
\end{equation}
From Step 1 we know that $f(x)$ tends to 0 as $|x|$ tends to infinity,
then the sequences $x_i$ and $y_i$, $i\in{\bf N}$, are bounded
and we may assume that there exist $x,y$ with $x,-y\in{\rm cl}(P)$, such
that 
\begin{equation*}
\label{II.4}
\lim_{i\to\infty}x_i=x\,,\;
\lim_{i\to\infty}y_i=y\,,
\end{equation*}
and obviously $z=\frac{1}{2}(x+y)$. On the other hand by (\ref{II.1b})
$$
f(x)\ge\limsup_{n\to\infty}f(x_i)\,,\;
f(-y)\ge\limsup_{n\to\infty}f(-y_i)\,.
$$
Thus
\begin{eqnarray*}
\Delta f(z)\ge\sqrt{f(x)f(-y)}\ge
\limsup_{n\to\infty}\sqrt{f(x_i)}\,
\limsup_{n\to\infty}\sqrt{f(-y_i)}\ge
\lim_{n\to\infty}\sqrt{f(x_i)f(-y_i)}=\Delta f(z)\,.
\end{eqnarray*}

\bigskip

{\em Step 4.} The function
\begin{equation*}
F(z)=\int_{{\bf R}^n}f(y)f(y-z)\,dy\,,\;y\in{\bf R}^n
\end{equation*}
is continuous in ${\bf R}^n$. 

{\em Proof.} Let $z, z_i\in{\bf R}^n$, $i\in{\bf N}$, be
such that $\lim_{i\to\infty}z_i=z$. Notice that, by the log-concavity, $f$ is 
continuous in ${\bf R}^n\setminus\partial P$, i.e. it is continuous a.e. in 
${\bf R}^n$. Thus the sequence of functions $f(y)f(y-z_i)$, $i\in{\bf N}$, 
converges to $f(y)f(y-z)$ for a.e. $y\in{\bf R}^n$. Moreover
$$
f(y)f(y-z_i)\le(\sup_{{\bf R}^n}f)\,f(y)\,,\quad
\forall x\in{\bf R}^n\;\mbox{and}\;\forall i\in{\bf N}\,.
$$
As $f\in L^1({\bf R}^n)$ and $\sup_{{\bf R}^n}f<\infty$ (Step 1) we may
apply the Dominated Convergence Theorem to obtain
$$
F(z)=\int_{{\bf R}^n}f(y)f(y-z)\,dy=\lim_{i\to\infty}
\int_{{\bf R}^n}f(y)f(y-z_i)\,dy=\lim_{i\to\infty} F(z_i)\,,
$$
i.e. $F$ is continuous at $z$.

\bigskip

{\em Step 5.} We prove that ${\rm int}(P)$ is a (convex) cone. 

{\em Proof.} As $f$ renders inequality (\ref{I.4}) en equality, we know from the
proof of Theorem \ref{teo.I.1}, and in particular from (\ref{II.0}),
that
\begin{equation}
\label{II.5}
F(z)=\frac{1}{2^n}\left(\int_{{\bf R}^n}f(y)\,dy\right)\Delta f(z)\,,
\end{equation}
for a.e. $z\in{\bf R}^n$. The function $\Delta f$ is strictly positive in the set
$$
\Delta P=\frac{1}{2}(P+(-P))=\left\{\frac{1}{2}(x+y)\;:\;x,-y\in P\right\}\,,
$$
and from ${\rm int}(P)\ne\emptyset$ it follows that $0\in{\rm int}(\Delta P)$, 
so that $\Delta f$ is continuous in a neighborhood of $0$ (recall that $\Delta f$ is 
log-concave). By this fact and Step 4, we deduce that equality 
(\ref{II.5}) is satisfied pointwise in a neighborhood of the origin.

From Step 1 and from Step 2 we easily get
$\sup_{{\bf R}^n}f=\max_{{\bf R}^n}f$. Up to a translation of the $x$
variable and to a multiplication of $f$ by a positive constant, we may
assume that
\begin{equation}
\label{II.5b}
\max_{x\in{\bf R}^n}f(x)=f(0)=1\,.
\end{equation}
A direct consequence of the previous equalities and of the definition of
$\Delta f$ is 
$$
\Delta f(0)=f(0)=1\,.
$$
Now, let us write equality (\ref{II.5}) for $z=0$:
\begin{equation*}
F(0)=\int_{{\bf R}^n}f^2(y)\,dy=
\frac{1}{2^n}\int_{{\bf R}^n}f(y)\,dy=
\int_{{\bf R}^n}f(2y)\,dy\,.
\end{equation*}
On the other hand, as $f$ is log-concave 
\begin{equation*}
f(y)\ge\sqrt{f(0)f(2y)}=\sqrt{f(2y)}\quad\mbox{for every $y\in{\bf R}^n$.}
\end{equation*}
We deduce that
\begin{equation}
\label{II.6}
f(y)=\sqrt{f(2y)}\quad\mbox{for a.e. $y\in{\bf R}^n$.}
\end{equation}
As a consequence of this fact, the sets $P$ and $2P$ coincide up to a 
null set (i.e. $\chi_P=\chi_{2P}$ a.e. in ${\bf R}^n$), but, since these sets are 
convex, this implies ${\rm int}(P)={\rm int}(2P)$, which is possible if and only 
if ${\rm int}(P)$ is a convex cone, with vertex at the origin. This concludes the 
proof of the assert of the present step. 

In the rest of the proof we will always assume that (\ref{II.5b}) holds
and that the vertex of $P$ is at the origin.


\bigskip

{\em Step 6.} We prove that the function $\Delta f$ is positive and
continuous in ${\bf  R}^n$. 

{\em Proof.} $\Delta f$ is positive in the set
\begin{equation*}
\frac{1}{2}(P+(-P))\supset\frac{1}{2}({\rm int}(P)+{\rm int}(-P))\,;
\end{equation*}
on the other hand, by the previous step, ${\rm int}(P)$ is a convex cone
with interior points, and an easy consequence of this fact is that
${\rm int}(P)+{\rm int}(-P)={\bf R}^n$. Thus $\Delta f$ is positive in
${\bf R}^n$ and, as it is log-concave, it is also continuous in ${\bf R}^n$.

\bigskip

{\em Step 7.} Let us fix an arbitrary $z\in{\bf R}^n$ and let $x_1$, $x_2$ be such 
that
\begin{equation*}
\frac{1}{2}(x_1+x_2)=z\,,\quad \Delta f(z)=\sqrt{f(x_1)f(-x_2)}
\end{equation*}
(see Step 3); notice that, as $\Delta f(z)>0$, $x_1$ and $-x_2\in P$. 
We prove that 
\begin{eqnarray}
\label{II.8}
&&f\left(\frac{1}{2}(2y-x_1)+\frac{1}{2}x_1\right)
f\left(\frac{1}{2}(2y-x_1)+\frac{1}{2}(-x_2)\right)=\nonumber\\
&&\\
&&\sqrt{f(2y-x_1)f(x_1)}\sqrt{f(2y-x_1)f(-x_2)}\,,\;
\mbox{for a.e. $y\in{\bf R}^n$.}\nonumber
\end{eqnarray}

{\em Proof.} The starting point is equation (\ref{II.5}); first notice that such
equation holds for every $z\in{\bf R}^n$, because 
it holds a.e. in ${\bf R}^n$ and $F(z)$ and $\Delta f(z)$ are
continuous in ${\bf R}^n$. On the other hand, (\ref{II.5}) can
be written as follows 
\begin{eqnarray*}
&&\int_{{\bf R}^n}f\left(\frac{1}{2}(2y-x_1)+\frac{1}{2}x_1\right)
f\left(\frac{1}{2}(2y-x_1)+\frac{1}{2}(-x_2)\right)\,dy\\
\\
&&=\int_{{\bf R}^n}\sqrt{f(2y-x_1)f(x_1)}\sqrt{f(2y-x_1)f(-x_2)}\,dy\,,
\end{eqnarray*}
so that equality (\ref{II.8}) follows immediately from the log-concavity of $f$.

\bigskip

{\em Step 8.} We prove the following fact: let $z$, $x_1$ and $x_2$ be as in the 
previous step, then for every point $x\in{\rm int}(P)$ the function $v=-\log f$ 
restricted either to the segment joining $x$ and $x_1$ or to the segment joining 
$x$ and $-x_2$, is affine, i.e.
\begin{eqnarray}
\label{II.10}
&&v((1-t)x+tx_1)=(1-t)v(x)+tv(x_1)\,,\;\forall t\in[0,1]\,,\nonumber\\
\\
&&v((1-t)x+t(-x_2))=(1-t)v(x)+tv(-x_2)\,,\;\forall t\in[0,1]\,.\nonumber
\end{eqnarray}

{\em Proof.} We write equation (\ref{II.8}) with $x=2y-x_1$:
\begin{eqnarray*}
f\left(\frac{1}{2}(x+x_1)\right)
f\left(\frac{1}{2}(x-x_2)\right)=
\sqrt{f(x)f(x_1)}\sqrt{f(x)f(-x_2)}\,,\;
\mbox{for a.e. $x\in{\bf R}^n$.}
\end{eqnarray*}
As $x_1$ and $-x_2$ belong to $P$, and as ${\rm int}(P)$ is a cone, the left 
hand-side and the right hand-side, as functions of $x$, are continuous in 
${\rm int}(P)$, so that the above equality holds pointwise for $x\in{\rm int}(P)$. 
Moreover, for $x\in{\rm int}(P)$ both sides of the equality are positive; on the
other hand we have, by the log-concavity of $f$,
\begin{eqnarray*}
f\left(\frac{1}{2}(x+x_1)\right)\ge
\sqrt{f(x)f(x_1)}\,,\quad
f\left(\frac{1}{2}(x-x_2)\right)\ge
\sqrt{f(x)f(-x_2)}\,.
\end{eqnarray*}
Thus, the last two inequalities must be in fact equalities; in terms of
$v$ this means that:
\begin{eqnarray*}
v\left(\frac{1}{2}(x+x_1)\right)=
\frac{1}{2}(v(x)+v(x_1))\,,\quad
v\left(\frac{1}{2}(x-x_2)\right)=
\frac{1}{2}(v(x)+v(-x_2))\,.
\end{eqnarray*}
As $v$ is convex, this concludes the proof of (\ref{II.10}).
In the rest of the proof we will denote the function $-\log f$ by $v$.

\bigskip

{\em Step 9.} We prove that for every $y\in{\bf R}^n$ and for every $t\ge 0$
$$
v(ty)=tv(y)\,.
$$

{\em Proof.} If $y\in{\bf R}^n\setminus{\rm cl}(P)$ the claim is obvious, since
$v(ty)=tv(y)=\infty$. If $y\in{\rm int}(P)$, the claim follows from the previous 
Step 8, where we take $z=0$ and consequently $x_1=x_2=0$. 
As last case, let $y\in\partial P$; by Step 2, and in particular
(\ref{II.1b}), there exists a sequence $y_i$,
$i\in{\bf N}$, contained in ${\rm int}(P)$, such that
$$
\lim_{i\to\infty}y_i=y\,,\quad\lim_{i\to\infty}v(y_i)=
\liminf_{x\to y,x\in{\rm int}(P)}v(x)=v(y)\,.
$$
Then clearly we also have
$$
tv(y)=
t\lim_{i\to\infty}v(y_i)=
t\liminf_{x\to y,x\in{\rm int}(P)}v(x)=
\liminf_{x\to y,x\in{\rm int}(P)}v(tx)=
\liminf_{x\to ty,x\in{\rm int}(P)}v(x)
=v(ty)\,,
$$
where we have used again the claims of Step 8 and Step 2.

\bigskip

{\em Step 10.} Let $z\in{\bf R}^n$ and $x_1$, $x_2$ be such that
\begin{equation*}
\frac{1}{2}(x_1+x_2)=z\quad{\rm and}\quad f(z)=\sqrt{f(x_1)f(-x_2)}\,.
\end{equation*}
We prove that
\begin{equation}
\label{II.9}
{\rm int}(P)\cap({\rm int}(P)+z)={\rm int}(P)+\frac{1}{2}x_1\quad{\rm and}\quad
{\rm cl}(P)\cap({\rm cl}(P)+z)={\rm cl}(P)+\frac{1}{2}x_1\,.
\end{equation}

{\em Proof.} From equation (\ref{II.8}) we obtain
\begin{equation*}
f(y)f(y-z)=f(2y-x_1)\Delta f(z)\,,\,\mbox{for a.e. $y\in{\bf R}^n$.}
\end{equation*}
The right hand-side is positive if and only if $y\in P\cap(P+z)$, while the left 
hand-side is positive if and only if $z\in\frac{1}{2}(D+x_1)=D+\frac{1}{2}x_1$ 
(recall that $\Delta f$ is positive in ${\bf R}^n$). Consequently
$$
\chi_{P\cap(P+z)}=\chi_{P+\frac{1}{2}x_1}\quad\mbox{a.e. in ${\bf R}^n$}
$$
from which equalities (\ref{II.9}) follow easily. 

\bigskip

{\em Step 11.} Let $z\in{\rm cl}(P)$ and $x_1, x_2$ as in the previous
step; we prove that  
$$
x_1=2z\quad{\rm and}\quad x_2=0\,.
$$

{\em Proof.} Using (\ref{II.9}) (second equality) we obtain that 
$z=0+z\in{\rm cl}(P)+\frac{1}{2}x_1$; as $z-\frac{1}{2}x_1=\frac{1}{2}x_2$,
this implies that $\frac{1}{2}x_2$ belongs to ${\rm cl}(P)$ and then, as 
${\rm cl}(P)$ is a cone, the same is true for $x_2$. On the other hand
$-x_2$ belongs to ${\rm cl}(P)$ also. Assume, by contradiction, that
$x_2\ne 0$, and let $r$ be the straight line through $x_2$ and $-x_2$. As
${\rm cl}(P)$ is a cone, it contains $r$. 
This implies in particular that if $x\in r$ and $w\in{\rm cl}(P)$,
then $x+s(w-x)\in{\rm cl}(P)$ for every $s\ge 0$ (indeed $(1-s)x\in
r\subset{\rm cl}(P)$ for every $s\in{\bf R}$ and $sw\in{\rm cl}(P)$ for
every $s\ge0$). Furthermore, if $w\in{\rm int}(P)$,
then $x+s(w-x)\in{\rm int}(P)$ for every $s>0$.

Let $y_i$, $i\in{\bf N}$, be a sequence of points contained in ${\rm
  int}(P)$ such that 
\begin{equation*}
\lim_{i\to\infty}y_i=0\,,\quad
\lim_{i\to\infty}v(y_i)=v(0)=0\,.
\end{equation*}
For every $i\in{\bf N}$, let $t_i=2y_i+x_2$; by the above considerations,
$t_i\in{\rm int}(P)$. 
By Step 8 we have, for every $i\in{\bf N}$,
\begin{equation*}
v(y_i)=v\left(\frac{1}{2}(-x_2+t_i)\right)=\frac{1}{2}(v(-x_2)+v(t_i))\,,
\end{equation*}
and, passing to the limit we get
\begin{equation*}
v(-x_2)+\lim_{i\to\infty}v(t_i)=0\,.
\end{equation*}
Since $t_i$ tends to $x_2$, by Step 2 and the previous equality we get 
$$
0\le v(-x_2)+v(x_2)\le v(-x_2)+\lim_{i\to\infty}v(t_i)=0
$$ 
(we recall that $v$ is nonnegative by (\ref{II.5b}). We deduce
$v(x_2)=v(-x_2)=0$ and this implies, by Step 9, that $v$ vanishes on the
whole line $r$. Thus $f(x)=1$ for every $x\in r$ 
and this contradicts the fact that the set $\{x\,|\,f(x)\ge1\}$ is bounded
established in Step 1. We infer that $x_2=0$ which concludes the proof
of the claim of the present step. 

\bigskip

{\em Step 12.} We show that $v$ is linear in ${\rm cl}(P)$. 
 
{\em Proof.} Take arbitrary points $x$ and $y$ in ${\rm int}(P)$, and
let $z=2x$; then 
by Step 11 $x$ coincides with the point $x_1$ corresponding to $z$ and
by Step 8 $v$ is affine on the segment joining $x$ and $y$. Then $v$ is
affine in ${\rm int}(P)$ and since $v(0)=0$, $v$ is linear in ${\rm
  int}(P)$. Moreover, for every $x\in\partial P$, by Step 2 we have 
$$
v(x)=\liminf_{y\to x,y\in{\rm int}(P)}v(y)
=\lim_{y\to x,y\in{\rm int}(P)}v(y)\,,
$$
so that $v$ is linear in ${\rm cl}(P)$.
\bigskip

{\em Step 13.} We prove that ${\rm cl}(P)$ does not contain any straight
line. 

{\em Proof.} By contradiction, let $r$
be a straight line contained in ${\rm cl}(P)$; as $v$ is linear and nonnegative
in ${\rm cl}(P)$, $v$ must be constant on $r$: $v(x)=c\ge 0$ for 
every $x\in r$. Consequently $f(x)\ge e^{-c}>0$ for every $x\in r$ which
contradicts claim {\em i}) of Step 1.

\bigskip

{\em Step 14.} We prove that the intersection of $P$ with a suitable
hyperplane is a $(n-1)$-dimensional simplex.

{\em Proof.} As ${\rm cl}(P)$ does not
contain any straight line, we may assume after a change of coordinates that
$$
{\rm cl}(P)\subset\{x=(x_1,\dots,x_n)\in{\bf R}^n\,:\,x_n>0\}\cup\{0\}\,.
$$
In particular, the section of ${\rm cl}(P)$:
$$
S={\rm cl}(P)\cap\{x=(x_1,\dots,x_n)\in{\bf R}^n\,:\,x_n=1\}
$$
is a compact convex set. We set $\pi=\{x=(x_1,\dots,x_n)\in{\bf
  R}^n\,:\,x_n=1\}$; let $(z_1,\dots,z_{n-1})$ be an arbitrary point
of ${\bf R}^{n-1}$ and let ${\bar z}=(z_1,\dots,z_{n-1},0)$. From Step 10 it 
follows that there exists $x\in{\bf R}^n$ such that
$$
{\rm cl}(P)\cap({\rm cl}(P)+{\bar z})={\rm cl}(P)+x\,;
$$
if we take the intersections with $\pi$ of both sides of the previous equality
we obtain
$$
S\cap(S+{\bar z})=({\rm cl}(P)+x)\cap\pi\,.
$$
Notice that, as ${\rm cl}(P)$ is a cone, the set $({\rm cl}(P)+x)\cap\pi$ is 
either empty or homothetic (i.e. equal up to a translation and a dilatation)
to $S$. Thus, we have proved that for any translate $S'$ of $S$
contained in $\pi$, $S\cap S'$ is either empty or is homothetic to
$S$. From Lemma 4 in \cite{Rogers-Shephard1} (see also \S 7.3 in \cite
{Schneider}) we obtain that $S$ is an $(n-1)$-dimensional simplex.

\bigskip

{\em Step 15.} Summarizing the conclusions of the previous steps, we
have proved that, up to a translation of the $x$ variable and a
multiplication of $f$ by a positive constant:

\begin{enumerate}

\item ${\rm int}(P_f)$ is an infinite convex cone with vertex at the
  origin having an $(n-1)$-dimensional simplex as a section;

\item $v=-\log f$ is linear in ${\rm int}(P_f)$.
\end{enumerate}

Claims {\em i}) and {\em ii}) are easy consequences of these facts.

\begin{flushright}
$\Box$
\end{flushright}

\section{An application to convex bodies}
\label{Paragrafo3}

In this section we give a new proof, based on (\ref{I.4}), of an
inequality proved by Rogers and Shephard in \cite{Rogers-Shephard2},
concerning the volume of the convex hull of the union of a convex body and its
reflected body with respect to the origin. For an
arbitrary set $A$, ${\rm co}(A)$ denotes the convex hull of $A$. 

\begin{thm}\label{teo.III.1}
{\bf (Rogers-Shephard.)} Let $K$ be a convex body in ${\bf R}^n$, then for
every $x\in K$ we have the following inequality 
\begin{equation}
\label{III.1}
V_n({\rm co}(K\cup(x-K)))\le2^nV_n(K)\,.
\end{equation}
\end{thm}

\bigskip

Inequality (\ref{III.1}) is optimal, indeed equality holds when $K$ is a
simplex and $x$ is one of its vertices; in \cite{Rogers-Shephard2} it is
proved that this is the only possibility. 

\bigskip

\noindent{\bf Proof of Theorem \ref{teo.III.1}.} We may assume that $0\in{\rm
  int}(K)\ne\emptyset$. Let us denote by $h_K$ the support function of $K$ and
by $K^*$ its polar with respect to the origin (for the definition of these notions we refer
to \cite{Schneider}); then $h_K$ is positive in ${\bf R}^n$ and $K^*$ is also a convex body.  

The following formula is a rather simple consequence of the fact that the
radial function of $K^*$ is the reciprocal of $h_K$ (see \cite{Schneider}, \S
1.7):  
\begin{equation}
\label{II.11}
\int_{{\bf R}^n}e^{-h_K(x)}\,dx=n!\,V_n(K^*)\,.
\end{equation}
The function $f=e^{-h_K}$ is log-concave as $h_K$ is convex (see \cite{Schneider}). We prove that
\begin{equation}
\label{II.12}
\Delta f=e^{-h_{K\cap(-K)}}\,.
\end{equation}
Indeed $\Delta f=e^{-\delta h_K}$, where $\delta h_K$ is defined by formula
(\ref{I.0aa}) in \S 1; $\delta h_K$ can be written in the form 
\begin{equation}
\label{II.13}
\delta h_K=\left[\frac{(h_K(\cdot))^*+(h_K(-\cdot))^*}{2}\right]^*
=\left[\frac{h_K^*+h_{-K}^*}{2}\right]^*\,,
\end{equation}
where ${}^*$ denotes the Legendre conjugation of convex functions (see Theorem
16.4 in \cite{Rockafellar}). On the other hand, the conjugate of the support
function of a convex body $L$ is the so-called indicatrix function ${\bf I}_L$
of $L$:     
\begin{equation*}
h_L^*(z)={\bf I}_L(z)=
\left\{
\begin{array}{ll}
\infty\quad&\mbox{if $z\in L$,}\\
0     \quad&\mbox{otherwise,}
\end{array}
\right.
\end{equation*}
and ${\bf I}_L^*=h_L$. Formula (\ref{II.12}) is then a consequence of
(\ref{II.13}). Using Theorem \ref{teo.I.1}, (\ref{II.11}) 
and (\ref{II.12}) we obtain
\begin{equation}
\label{II.14}
V_n((K\cap(-K))^*)\le2^nV_n(K^*)\,.
\end{equation}
Moreover $(K\cap(-K))^*={\rm co}(K^*\cup(-K^*))$, by Theorem 1.6.2 in
\cite{Schneider}. Thus (\ref{II.14}) becomes $V_n({\rm
  co}(K^*\cup(-K^*)))\le2^nV_n(K^*)$ and as $K$ is arbitrary we deduce  
\begin{equation*}
V_n({\rm co}(K\cup(-K)))\le2^nV_n(K)
\end{equation*}
for every convex body $K$ such that $0\in{\rm int}(K)$. Clearly, we could have chosen any interior point of $K$
as origin so that (\ref{III.1}) is valid for every $x$ in the interior of $K$;
in the general case the formula follows by a continuity argument. 
\begin{flushright}
$\Box$
\end{flushright}


\section{Difference functions of order $\alpha$ and related inequalities}
\label{Paragrafo4}

Throughout this section, $\alpha$ is a parameter varying in $[-\infty,0]$. We
start by defining the mean of order $\alpha$ of two nonnegative numbers. For
$a,b\ge0$ and $t\in[0,1]$, we set 
\begin{equation}
\label{IV.1}
M_\alpha(a,b;t)=
\left\{
\begin{array}{llll}
a^t\,b^{1-t}&\;\mbox{if}\;\alpha=0\,,\\
\left(ta^\alpha+(1-t)b^\alpha\right)^{1/\alpha}&\;\mbox{if}\;\alpha\in(-\infty,0),\;a\ne0,\;b\ne0\,,\\
0&\;\mbox{if}\;\alpha\in(-\infty,0)\;\mbox{and either $a=0$ or $b=0$}\,,\\
\min\{a,b\}&\;\mbox{if}\;\alpha=-\infty\,.
\end{array}
\right.
\end{equation}
For $\alpha<0$ this definition can be extended naturally to the case in which either $a=\infty$ or
  $b=\infty$, setting $M_\alpha(\infty,c)=M_\alpha(c,\infty)=\displaystyle{\frac{c}{2^{1/\alpha}}}$
  for every $c\ge0$, and $M_\alpha(\infty,\infty)=\infty$. 

\begin{de} We say that a function $f\,:{\bf R}^n\rightarrow[0,\infty]$ is $\alpha$-concave if  
\begin{equation*}
f(tx+(1-t)y)\ge M_\alpha(f(x),f(y);t)\,,\quad\forall x,y\in{\bf R}^n\,,\;\forall t\in[0,1]\,.
\end{equation*}
\end{de}

For $\alpha=0$ we get log-concave functions. 
For $\alpha=-\infty$,
$\alpha$-concave functions are called {\em quasi-concave} functions; let us point out
the following well-known characterization. 

\begin{prop}
\label{prop.IV.1}
A function $f\,:{\bf R}^n\rightarrow[0,\infty]$ is quasi-concave if and only
if for every $s\ge 0$ the set $\{x\in{\bf R}^n\,:\,f(x)> s\}$ is convex. 
\end{prop}

\bigskip

We define the {\em difference function of order $\alpha$} of a function
$f\,:{\bf R}^n\rightarrow[0,\infty]$ in the following way 
\begin{equation}
\label{IV.2}
\Delta_\alpha
f(z)=\sup\left\{M_\alpha\left(f(x),f(-y);\frac{1}{2}\right)\,:\,\frac{1}{2}(x+y)=z\right\}\,,\quad
z\in{\bf R}^n\,. 
\end{equation}
For $\alpha=0$ we retrieve the definition of difference function given in \S 1. 

\begin{prop}
\label{prop.IV.2}
If $f$ is $\alpha$-concave, then $\Delta_\alpha f$ is also $\alpha$-concave. 
\end{prop}
{\bf Proof.} The case $\alpha=0$ has already been considered in \S 1 and the
argument for $\alpha\in(-\infty,0)$ is very similar; indeed in this case we
have 
\begin{equation*}
(\Delta_\alpha
f)^\alpha(z)=\sup\left\{\frac{f^\alpha(x)+f^\alpha(-y)}{2}\,:\,\frac{1}{2}(x+y)=z\right\}\,,\quad
z\in{\bf R}^n\,, 
\end{equation*}
i.e. $(\Delta_\alpha f)^\alpha$ is the infimal convolution of
$f^\alpha(\cdot)$ and $f^\alpha(-\cdot)$, which are convex (as $f$ is
$\alpha$-concave and $\alpha<0$); hence $(\Delta_\alpha f)^\alpha$ is also
convex (again, we refer to \cite{Rockafellar}, Chapter 5) and consequently
$\Delta_\alpha f$ is $\alpha$-concave. 

In the case $\alpha=-\infty$ we start from the following equality, which will
be helpful also in the sequel: for every $s\ge0$
\begin{equation}
\label{IV.3}
\{z\in{\bf R}^n\,:\,\Delta_{-\infty}f(z)> s\}=
\frac{1}{2}\{x\in{\bf R}^n\,:\,f(x)> s\}+
\frac{1}{2}\{y\in{\bf R}^n\,:\,f(-y)> s\}\,.
\end{equation}   
In order to prove it, assume that $z\in\{z\in{\bf
  R}^n\,:\,\Delta_{-\infty}f(z)> s\}$; then there exist $x$ and $y$
such that $\min\{f(x),f(-y)\}>s$ and $\dfrac{1}{2}(x+y)=z$,
i.e. $z\in\dfrac{1}{2}\{x\in{\bf
  R}^n\,:\,f(x)>s\}+\dfrac{1}{2}\{y\in{\bf
  R}^n\,:\,f(-y)>s\}$. Then
\begin{eqnarray*}
\{z\in{\bf R}^n\,:\,\Delta_{-\infty}f(z)> s\}\subset
\frac{1}{2}\{x\in{\bf R}^n\,:\,f(x)> s\}+\frac{1}{2}\{y\in{\bf R}^n\,:\,f(-y)> s\}.
\end{eqnarray*}
The reverse inclusion can be proved in a similar way. By Proposition
\ref{prop.IV.1} the two sets on the left hand-side of (\ref{IV.3}) are convex
and then the set on the right hand-side is also convex; hence
$\Delta_{-\infty}f$ is quasi concave again by Proposition \ref{prop.IV.1}.
\begin{flushright}
$\Box$
\end{flushright}

\bigskip

\noindent{\bf Problem.} {\em For $\alpha\le0$ and $n\in{\bf N}$, determine the number}
\begin{equation*}
C(n,\alpha)=\sup\left\{\dfrac{\displaystyle{\int_{{\bf R}^n}\Delta_\alpha
      f\,dz}}{\displaystyle{\int_{{\bf R}^n}f\,dx}}\,|\,f\,:\,{\bf
    R}^n\rightarrow[0,\infty]\,,\;f\;\mbox{$\alpha$-concave}\,,\;0<
\displaystyle{\int_{{\bf R}^n}f\,dx}<\infty\right\}\,.    
\end{equation*}

\bigskip

Theorems \ref{teo.I.1} and \ref{teo.I.2} provide the solution
for $\alpha=0$: $C(n,0)=2^n$ for every $n$. In the remaining part of the
paper we will provide some further partial answers to this problem. 
Let us start from the case $\alpha=-\infty$. The following theorem is a
relatively simple application of the Rogers-Shephard inequality.

\begin{thm}\label{teo.IV.1} Let $f$ be a quasi-concave function such that
  $f\in L^1({\bf R}^n)$. Then 
\begin{equation}
\label{IV.4}
\int_{{\bf R}^n}\Delta_{-\infty}f(x)\,dx\le
\frac{1}{2^n}\binom{2n}{n} \int_{{\bf R}^n}f(x)\,dx\,.
\end{equation}
\end{thm}
\noindent{\bf Proof.} For every $s>0$ the set 
$\{z\in{\bf R}^n\,:\,f(z)>s\}$ is convex (by Proposition \ref{prop.IV.1}) and it
is bounded, by the assumption $f\in L^1({\bf R}^n)$. Then  
$\{z\in{\bf R}^n\,:\,\Delta_{-\infty}f(z)>s\}$ is also convex and, by
(\ref{IV.3}), bounded. Using (\ref{IV.3}) and the Rogers-Shephard inequality we
obtain
$$
V_n(\{z\in{\bf R}^n\,:\,\Delta_{-\infty}f(z)> s\})\le
\frac{1}{2^n}\binom{2n}{n}
V_n(\{x\in{\bf R}^n\,:\,f(x)> s\})\,,\quad\forall s>0\,.
$$
Using the layer cake principle and the previous inequality we obtain (\ref{IV.4}).
\begin{flushright}
$\Box$
\end{flushright}

Inequality (\ref{IV.4}) is sharp; indeed it is clear from its proof and from
equality cases in the Rogers-Shephard inequality, that if $f$ is such that 
$\{x\in{\bf R}^n\,:\,f(x)>s\}$ is either empty or a simplex
for every $s\ge0$, then equality holds in (\ref{IV.4}) (we can take, for
instance, $f$ to be the characteristic function of a simplex). Then we can
state the following

\begin{cor} For every $n\in{\bf N}$, $C(n,-\infty)=
\displaystyle{\frac{1}{2^n}\binom{2n}{n}}$.
\end{cor}

The case $\alpha=-\infty$ permits to find an upper bound for $C(n,\alpha)$ in
the general case.

\begin{thm} For every $\alpha<0$ and $n\in{\bf N}$, $C(n,\alpha)\le
\displaystyle{\frac{1}{2^{n+1/\alpha}}\binom{2n}{n}}$.
\end{thm}

\noindent{\bf Proof.} Let $f$ be a nonnegative $\alpha$-concave function such
that $f\in L^1({\bf R}^n)$. By the monotonicity property of the mean of order
$\alpha$ with respect to $\alpha$, $f$ is also quasi-concave. 
Moreover, it is quite an easy exercise to prove that for every $a$, $b\ge0$ 
\begin{equation}\label{IV.5}
M_\alpha\left(a,b,\frac{1}{2}\right)\le
\frac{1}{2^{1/\alpha}}\min(a,b)\,.
\end{equation}
Then $\Delta_\alpha
f\le\displaystyle{\frac{1}{2^{1/\alpha}}}\Delta_{-\infty}f$ in ${\bf R}^n$;
using this fact and Theorem \ref{teo.IV.1} we obtain
\begin{equation*}
\int_{{\bf R}^n}\Delta_{\alpha}f(x)\,dx\le
\frac{1}{2^{1/\alpha}}
\int_{{\bf R}^n}\Delta_{-\infty}f(x)\,dx\le
\frac{1}{2^{1/\alpha}}
\frac{1}{2^n}\binom{2n}{n} \int_{{\bf R}^n}f(x)\,dx\,.
\end{equation*}
\begin{flushright}
$\Box$
\end{flushright}

\section{The one dimensional case}
\label{Paragrafo5}

The constant $C(1,\alpha)$ can be determined explicitely for every $\alpha$.

\begin{thm}\label{teo.V.1} Let $\alpha\in(-\infty,0)$ and $f$ be an
  $\alpha$-concave function. Then
\begin{eqnarray}
\int_{\bf R}\Delta_\alpha f\,dx&\le&2\int_{\bf R}f\,dx\,,\quad{\rm if}\;\alpha\in(-1,0)\,,\label{V.1a}\\
\int_{\bf R}\Delta_\alpha
f\,dx&\le&\displaystyle{\frac{1}{2^{1/\alpha}}}\int_{\bf R}f\,dx\,,\quad{\rm
  if}\;\alpha\in(-\infty,-1]\,,\label{V.1b} 
\end{eqnarray}
Both inequalities are optimal, i.e. for every $\alpha\in(-1,0)$
(respectively, $\alpha\in(-\infty,-1]$) there exists an $\alpha$-concave function $f$
such that in (\ref{V.1a}) (respectively, in (\ref{V.1b})) equality holds. 
\end{thm}

As a consequence, in the notation introduced in the previous section,
\begin{equation*}
C(1,\alpha)=
\left\{
\begin{array}{ll}
2\quad&\mbox{for}\quad\alpha\in(-1,0)\,,\\
\displaystyle{\frac{1}{2^{1/\alpha}}}\quad&\mbox{for}\quad\alpha\in(-\infty,-1]\,.
\end{array}
\right.
\end{equation*}

The crucial ingredient in the proof of Theorem \ref{teo.V.1} is to prove that
for an arbitrary $f$ 
\begin{equation*}
\frac{\displaystyle{\int_{\bf R}\Delta_\alpha f\,dx}}{\displaystyle{\int_{\bf R}f\,dx}}\le
\frac{\displaystyle{\int_{\bf R}\Delta_\alpha f^*\,dx}}{\displaystyle{\int_{\bf R}f^*\,dx}}
\end{equation*}
where $f^*$ is a suitable decreasing rearrangement of $f$. 
Then we will prove inequalities (\ref{V.1a}) and (\ref{V.1b}) for
$f^*$. The latter task is particularly easy due to the features of the
function $f^*$. The definition and some properties of $f^*$ are contained in
the following lemma.

\begin{lm}\label{lemma.V.1} Let $\alpha<0$ and $f$ be an $\alpha$-concave
  function. Define
\begin{equation*}
f^*(z)=\left\{
\begin{array}{ll}
\sup_{x\in{\bf R}}\min\{f(x),f(x-z)\}\,,\quad&\mbox{if $z\ge0$}\,,\\
0                                            &\mbox{if $z<0$}\,.
\end{array}
\right.
\end{equation*}
Then
\begin{enumerate}

\item[i.] $f^*$ is $\alpha$-concave;
\item[ii.] for every $s\ge0$ 
$$
V_1(\{z\in{\bf R}\,:\,f^*(z)>s\})=
V_1(\{z\in{\bf R}\,:\,f(z)>s\})
$$
($V_1$ denotes the one-dimensional Lebesgue measure);
\item[iii.] 
$$
\int_{\bf R}f^*\,dz=
\int_{\bf R}f\,dz\,;
$$
\item[iv.] $f^*$ is decreasing in $[0,\infty)$;
\item[v.] $\Delta_\alpha
  f^*(z)=M_\alpha\left(f^*(0),f^*(2|z|),\displaystyle{\frac{1}{2}}\right)$ for
  every $z\in{\bf R}$;
\item[vi.] $\Delta_\alpha f^*(z)\ge\Delta_\alpha f(z)$ for every $z\in{\bf R}$.
\end{enumerate}
\end{lm}

\noindent{\bf Proof.} Proof of {\em i.} Let $z_1, z_2\in{\bf R}$ and
$t\in[0,1]$; we have to prove that
\begin{equation}\label{V.2}
f^*(tz_1+(1-t)z_2)\ge M_\alpha(f^*(z_1),f^*(z_2),t)\,.
\end{equation}
If either $z_1<0$ or $z_2<0$ this is true by the definition of $f^*$; hence,
let us assume that $z_1,z_2\ge0$. Let $x_1,x_2\in{\bf R}$;
\begin{eqnarray*}
f^*(tz_1+(1-t)z_2)&=&\sup_{x_1,x_2}\min\left\{f(tx_1+(1-t)x_2),f(t(x_1-z_1)+(1-t)(x_2-z_2))\right\}\\
                &\ge&\sup_{x_1,x_2}\min\left\{M_\alpha(f(x_1),f(x_2),t),M_\alpha(f(x_1-z_1),f(x_2-z_2),t)\right\}\\
                &\ge&\sup_{x_1,x_2}M_\alpha(\min\{f(x_1),f(x_1-z_1)\},\min\{f(x_2),f(x_2-z_2)\},t)\\
                  &=&M_\alpha(\sup_{x_1}\min\{f(x_1),f(x_1-z_1)\},\sup_{x_2}\min\{f(x_2),f(x_2-z_2)\},t)\\
                  &=&M_\alpha(f^*(z_1),f^*(z_2),t)\,,
\end{eqnarray*}
where we have used (in the first inequality) the $\alpha$-concavity of $f$.

Proof of {\em ii.} For $s\ge 0$ set 
$$
{\cal F}_s=\{z\in{\bf R}\,:\,f(z)>s\}\,,\;
{\cal F}^*_s=\{z\in{\bf R}\,:\,f^*(z)>s\}\,.
$$
Assume that $z\in{\cal F}^*_s$; then $z\ge 0$, moreover there exists $x$ such that
$\min\{f(x),f(x-z)\}>s$, i.e. $f(x)>s$ and $f(-(z-x))>s$ so that
$z=x+(z-x)\in({\cal F}_s+(-{\cal F}_s))\cap[0,\infty)$ and then 
${\cal F}^*_s\subset({\cal F}_s+(-{\cal F}_s))\cap[0,\infty)$. The reverse inclusion 
can be proved in a similar way. Then we have that 
\begin{equation}\label{V.2a}
{\cal F}_s^*=
\left({\cal F}_s+(-{\cal F}_s)\right)\cap[0,\infty)\,.
\end{equation}
The set ${\cal F}_s+(-{\cal F}_s)$ is symmetric with respect
to $z=0$ and its measure equals $2V_1({\cal F}_s)$; then (\ref{V.2a}) implies that
$V_1({\cal F}_s^*)=V_1({\cal F}_s)$ for every $s$.  

Proof of {\em iii.} This is an immediate consequence of {\em ii.} and the
layer cake principle.

Proof of {\em iv.} Let $0\le z_1\le z_2$, $x\in{\bf R}$. There exists $t\in[0,1]$
such that $x-z_1=t(x-z_2)+(1-t)x$. From the $\alpha$-concavity of $f$ it
follows that $f(x-z_1)\ge M_\alpha(f(x-z_2),f(x),t)\ge\min\{f(x),f(x-z_2)\}$. 
Hence $\min\{f(x),f(x-z_1)\}\ge\min\{f(x),f(x-z_2)\}$ for every $x$; this
yields $f^*(z_1)\ge f^*(z_2)$.

Proof of {\em v.} Let $z\ge 0$. 
\begin{eqnarray*}
\Delta_\alpha f^*(z)&=&\sup_{x\in{\bf R}}M_\alpha\left(f^*(x),f^*(x-2z),\frac{1}{2}\right)
                     =\sup_{x\ge2z}M_\alpha\left(f^*(x),f^*(x-2z),\frac{1}{2}\right)\\
                    &=&M_\alpha\left(f^*(0),f^*(2z),\frac{1}{2}\right)\,,
\end{eqnarray*}
where in the last equality we have used the monotonicity of $f^*$. As
$\Delta_\alpha f^*$ is an even function, the claim of follows.

Proof of {\em vi.} For $z\ge0$ and $x\in{\bf R}$ we have
\begin{eqnarray*}
&&f^*(0) =\sup_{\bf R}f\ge\max\{f(x),f(x-2z)\}\,;
\\
&&f^*(2z)=\sup_{\xi\in{\bf R}}\min\{f(\xi),f(\xi-2z)\}\ge\min\{f(x),f(x-2z)\}\,.
\end{eqnarray*}
Hence,
\begin{equation*}
M_\alpha\left(f(x),f(x-2z),\frac{1}{2}\right)\le
M_\alpha\left(f^*(0),f^*(2z),\frac{1}{2}\right)\,;
\end{equation*}
taking the supremum of the left hand-side with respect to $x$ we obtain
$\Delta_\alpha f(z)\le\Delta_\alpha f^*(z)$ for every $z\ge 0$; as the
functions are even, the inequality is valid for every $z\in{\bf R}$.
\begin{flushright}$\Box$
\end{flushright} 

\bigskip

\noindent{\bf Proof of Theorem \ref{teo.V.1}.} By the previous lemma, for an
arbitrary $\alpha$-concave function $f$ we have 
\begin{equation}\label{V.3}
\frac{\displaystyle{\int_{\bf R}\Delta_\alpha f\,dx}}{\displaystyle{\int_{\bf
      R}f\,dx}}\le
\frac{\displaystyle{\int_{\bf R}\Delta_\alpha f^*\,dx}}{\displaystyle{\int_{\bf
      R}f^*\,dx}}=
2\,\frac{\displaystyle{\int_0^\infty
    M_\alpha\left(f^*(0),f^*(2x),\frac{1}{2}\right)\,dx}}{\displaystyle{\int_0^\infty f^*\,dx}}\,. 
\end{equation}
{\em Case $-1<\alpha<0$.} As $f^*$ is $\alpha$-concave we have  
$$
M_\alpha\left(f^*(0),f^*(2x),\frac{1}{2}\right)\le f^*(z)\,,\quad\forall\,
z\ge0\,.
$$
Inequality (\ref{V.1a}) follows immediately from the previous inequality and
(\ref{V.3}). In order to prove that (\ref{V.1a}) may be an equality, define a
function $\tilde{f}$ as follows   
\begin{equation}
\tilde{f}(x)=
\left\{
\begin{array}{ll}
(1+x)^{1/\alpha}\,,\quad&\mbox{if $x\ge0$,}\\
0\,,\quad&\mbox{if $x<0$.}
\end{array}
\right.
\end{equation}
We have: $\tilde{f}\in L^1({\bf R})$, $\tilde{f}$ is $\alpha$-concave and, as it is monotone decreasing in
$[0,\infty)$ and it vanishes in $(-\infty,0)$, $\tilde{f}=\tilde{f}^*$ in
${\bf R}$. Moreover, by a direct computation we see that
$$
M_\alpha\left(\tilde{f}(0),\tilde{f}(2x),\frac{1}{2}\right)=\tilde{f}(x)\,,\quad\forall\, x\ge0\,;
$$
i.e. $\Delta_\alpha\tilde{f}(x)=\tilde{f}(x)$ for $x\ge0$. Consequently,
(\ref{V.1a}) becomes an equality if $f=\tilde{f}$.\\
{\em Case $\alpha\le-1$.} Inequality (\ref{V.1b}) follows directly from
(\ref{V.3}) and 
\begin{equation}\label{V.4}
M_\alpha\left(f^*(0),f^*(2x),\frac{1}{2}\right)\le 
M_\alpha\left(\infty,f^*(2x),\frac{1}{2}\right)=\frac{1}{2^{1/\alpha}}f^*(2x)\,.
\end{equation}
Let us construct a function $\tilde{f}$ for which (\ref{V.1b}) is an
equality. We define
$$
\tilde{f}(x)=
\left\{
\begin{array}{llll}
0\,,           \quad&\mbox{if $x\ge1$,}\\
x^{1/\alpha}\,,\quad&\mbox{if $x\in(0,1)$,}\\
\infty      \,,\quad&\mbox{if $x=0$,}\\
0\,,           \quad&\mbox{if $x<0$.}\\
\end{array}
\right.
$$
As in the previous case, we see that: $\tilde{f}\in L^1({\bf R})$, $\tilde{f}$
is $\alpha$-concave and $\tilde{f}=\tilde{f}^*$ in ${\bf R}$. Moreover, for $f=\tilde{f}$
the inequality in (\ref{V.4}) is in fact an equality, so that
\begin{equation*}
\Delta_\alpha\tilde{f}(x)=\frac{1}{2^{1/\alpha}}\tilde{f}(2x)\,,\;\,\forall
x\ge0\,.
\end{equation*}
Then, for $f=\tilde{f}$ (\ref{V.1b}) is an equality.
\begin{flushright}$\Box$
\end{flushright}

\end{document}